\documentclass[fleqn,preprint,3p,a4paper]{elsarticle}

\usepackage{amssymb}
\usepackage{amsmath}
\usepackage{amsthm}
\usepackage{dcolumn}
\usepackage{endnotes}
\usepackage{tabularx}
\usepackage[matrix,arrow]{xy}
\usepackage{wasysym}

\newtheorem{theorem}{Theorem}[section]
\newtheorem{proposition}[theorem]{Proposition}
\newtheorem{lemma}[theorem]{Lemma}
\newtheorem{corollary}[theorem]{Corollary}

\theoremstyle{definition}
\newtheorem{definition}[theorem]{Definition}
\newtheorem{example}[theorem]{Example}
\newtheorem{remark}[theorem]{Remark}

\newcommand{\ir}{{\mathsf{Irr}}}

\newcommand{\SI}{{\rm SI}}

\newcommand{\cl}{{\rm cl}}
\newcommand{\ii}{{\rm int}}

\newcommand{\ua}{\mathord{\uparrow}}
\newcommand{\da}{\mathord{\downarrow}}

\journal{Topology and its applications}

\begin{document}

\begin{frontmatter}



\title{Characterizing $T_0$-spaces where GIS-convergence is topological \tnoteref{t1}}
\tnotetext[t1]{This research was supported by the National Natural Science Foundation of China (Nos. 12471070, 12071199, 12261040).}


 \author[X. Wen]{Xinpeng Wen}
\ead{wenxinpeng2009@163.com}

\author[F. Zhang]{Fugang Zhang}
\ead{zfgzhangfugang@163.com}

\author[J. He]{Jiawei He\corref{mycorrespondingauthor}}
\cortext[mycorrespondingauthor]{Corresponding author.}
\ead{hjwywh@mails.ccnu.edu.cn}

\address{School of Mathematics and Information Science, Nanchang Hangkong University, Nanchang 330063, China}

\begin{abstract}
  In this paper, we introduce and investigate the GIS-convergence in $T_0$-spaces and the related concept of GIS-continuous spaces(GSI-continuous spaces). It is proved that $SI(X)$ is a locally hypercompact space iff  $(X, \tau)$ is a $GSI$-continuous space iff GSI-convergence in $X$ with respect to the topology $\mathcal{O}_{SI}(X)$ is topological. We also show that the $GIS$-convergence in $(X, \tau)$ is topological if and only if $(X,\tau)$ is a $GIS$-continuous space.

\end{abstract}




%% MSC codes here, in the form:

\end{frontmatter}


\section{Introduction}

Convergence and convergence class play an important role in both order theory and general topology (see \cite{GHKLMS, K}). For a topological space $(X, \tau)$ and a class $\mathcal L$ consisting of pairs $((x_i)_{i\in I}, x)$, where $(x_i)_{i\in I}$ is a net in $X$ and $x$ a point of $X$, the topology $\tau$ can naturally induce a convergence class as follows:

\begin{center}
$\mathcal{C}(\tau)=\{((x_i)_{i\in I}, x) : (x_i)_{i\in I} \mbox{ is a net in } X, x\in X \mbox{ and for any } U\in\mathcal \tau$, $x\in U$ implies that $(x_i)_{i\in I} \mbox{ is eventually in } U\}$.
\end{center}

\noindent And we can define a topology on $X$ associated with $\mathcal L$:

\begin{center}
$\mathcal O(\mathcal L)=\{U\subseteq X : ((x_i)_{i\in I}, x)\in\mathcal L$ and $x\in U$ imply
 $x_i\in U$ eventually$\}$.
\end{center}
It is easy to verify that $\tau=\mathcal O(\mathcal C(\tau))$. However, if $\mathcal L$ is not a convergence class in the sense of Kelley \cite{K}, then the convergence class $\mathcal C(\mathcal O(\mathcal L))\neq\mathcal L$, that is, the class $\mathcal L$ is not topological.

Various types of convergence have been extensively explored by many researchers (see \cite{ASHZ, CK, E1, E, GHKLMS, LZ, SLL, ZBX, ZLW, ZZ, ZL}). By utilizing these distinct notions of convergence, not only have they introduced corresponding concepts of continuity for posets (and more generally, for topological spaces), but they have also revealed connections between order theory and topology. A fundamental and well-known theorem states that the lim-inf convergence in a dcpo $P$ is topological if and only if $P$ is a continuous domain (see \cite[Theorem II-1.9]{GHKLMS}). This theorem was subsequently extended to general partially ordered sets in \cite{ZZ}. Furthermore, Chen and Kou \cite{CK} defined the GS-convergence in a poset $P$ and demonstrated that it is topological precisely when $P$ is a quasi-continuous domain.

During his invited talk at the Sixth International Symposium on Domain Theory in 2013, Jimmie Lawson pointed out the necessity of developing the core of domain theory directly within $T_0$-spaces rather than restricting it to posets. Following this approach, Zhao and Ho \cite{ZH} proposed a significant new topology on $T_0$-spaces, known as the SI-topology, by replacing directed sets with irreducible sets; this serves as a generalization of the Scott topology on posets. Subsequently, Andradi et al. \cite{ASHZ} introduced both the SI-convergence in $T_0$-spaces and the notion of I-continuous spaces. They demonstrated that for a $T_0$-space $X$, the SI-topology on $X$ is induced by the SI-convergence, and provided that $X$ satisfies a certain condition, the SI-convergence in $X$ is topological precisely when $X$ is an I-continuous space. In \cite{ZLL},
It define and study $IS$-convergence in $T_0$ spaces and it prove that the IS-convergence in a $T_0$ space is topological if and only if the $T_0$ space is IS-continuous.

In this paper, motivated by the theory of SI-topologies, IS-convergence and IS-continuous spaces established in \cite{ASHZ, CK, LZ, ZH}, we introduce and investigate the GIS-convergence in $T_0$-spaces and the related concept of GIS-continuous spaces(GSI-continuous spaces). Our main result is that for a  $T_0$-space $X$, $SI(X)$ is a locally hypercompact space iff  $(X, \tau)$ is a $GSI$-continuous space iff GSI-convergence in $X$ with respect to the topology $\mathcal{O}_{SI}(X)$ is topological. We also show that the $GIS$-convergence in $(X, \tau)$ is topological if and only if $(X,\tau)$ is a $GIS$-continuous space.

\section{Preliminaries}

In this section, we briefly recall some basic concepts and results about ordered structures and $T_0$-spaces that will be used in the paper. For further details, we refer the reader to \cite{RE, GHKLMS, GL, K}.

For a poset $P$ and $A \subseteq P$, define $\ua A=\{x\in P: a\leq x \mbox{ for some }a\in A \}$ and $\da A=\{x \in P: x \leq a \mbox{ for some }a\in A\}$. For $x\in X$, let $\ua x=\ua\{x\}$ and $\da x=\da\{x\}$. A subset $A$ is called a \emph{lower set} (resp., an \emph{upper set}) if $A=\da A$ (resp., $A=\ua A$). If the set of upper bounds of~$A$ has a unique smallest element (that is, the set of upper bounds contains exactly one of its lower bounds), we call this element the \emph{least upper bound} and write it as~$\bigvee A$ or sup~$A$ (for supremum). Similarly the greatest lower bound is written as~$\bigwedge A$ or inf~$A$ (for infimum). A nonempty subset $D$ of a poset $P$ is called \emph{directed} if every finite subset of $D$ has an upper bound in $D$. The set of all directed sets of $P$ is denoted by $\mathcal{D}(P)$. The poset $P$ is called a \emph{directed complete poset}, or \emph{dcpo} for short, if for any
$D\in \mathcal D(P)$, $\bigvee D$ exists in $P$. For a nonempty subset $C$ of $P$, define $max(C)=\{c\in C : c \mbox{~ is a maximal element of~} C\}$ and $min(C)=\{c\in A : c \mbox{~ is a minimal element of~} C\}$.

For two subsets $F, G$ of a poset $P$, we say $F$ approximates $G$ if for every directed subset $D\subseteq P$, whenever $\bigvee D$ exists, $\bigvee D\in \uparrow G$ implies $d\in \uparrow F$ for some $d\in D$. A dcpo $P$ is called a \emph{quasi}-\emph{continuous domain} if for all $x\in P$, $\uparrow x$ is the filtered (with respect to set inclusion order) intersection of sets of the form $\uparrow F$, where $F$ approximates $\{x\}$ and $F$ is a nonempty finite set of $P$.

A subset $U$ of a poset $P$ is \emph{Scott open} if (i) $U=\mathord{\uparrow}U$, and (ii) for any directed subset $D$ for which $\bigvee D$ exists, $\bigvee D\in U$ implies $D\cap U\neq\emptyset$. The topology formed by all the Scott open sets of $P$ is called the \emph{Scott topology}, written as $\sigma(P)$. The upper sets of $P$ form the (\emph{upper}) \emph{Alexandroff topology} $\alpha(P)$. The lower topology on a poset $P$, generated by $\{P\setminus\uparrow x:x\in P\}$ (as a subbase), is denoted by $\omega(P)$. Dually, define the upper topology on $P$ and denote it by $v(P)$.

For a $T_0$-space $X$ and $A\subseteq X$, let $\mathcal{O}(X)$ (resp., $\Gamma(X)$) be the set of all open subsets (resp., closed subsets) of $X$. The closure of a subset $A$ in $X$ is denoted by $\cl_{X}~A$ or simply by $\cl~\!A$ and the interior of $A$ in $X$ by $\ii_{X} A$ or simply by $\ii ~ A$. A nonempty subset $F$ of $X$ is called an \emph{irreducible} set if for any $C_1, C_2\in\Gamma(X)$, $F\subseteq C_1\cup C_2$ implies $F\subseteq C_1$ or $F\subseteq C_2$. We denote by $\ir(X)$  the set of all irreducible  subsets of $X$. The space $X$ is is said to be \emph{irreducible}-\emph{complete} if every irreducible subset of $X$ has a sup. We use $\leq_X$ to denote the \emph{specialization order} of $X$: $x\leq y$ iff $x\in \cl~\! \{y\} $.
A space $X$ is called a $d$-space (or monotone convergence space) if $X$ (with the specialization order) is a dcpo and $\mathcal{O}(X)\subseteq \sigma(X)$.
In the following, when a $T_0$-space is considered as a poset, the order always refers to the specialization order if no other explanation is given. Clearly, all open sets  (resp., closed sets) of $X$ are upper sets (resp., lower sets) of $X$. Let $\mathbf{up}(X)=\{\uparrow_X U : U\subseteq X\}$ (that is, the family of all upper subsets of $X$), $X^{(<\omega)}$ be the set of all nonempty finite subsets of $X$ and $\mathbf{Fin}~\!\!X=\{\ua F : F\in X^{(<\omega)}\}$. For a subset $M$ of $X$, the \emph{induced topology} or \emph{subspace topology} on $M$ is denoted by $\mathcal O(X)|_M$, that is, $\mathcal O(X)|_M=\{U\cap M : U\in \mathcal O(X)\}$. A nonempty family $\mathcal{G}$ of subsets of a set $X$ is called \emph{filtered} if given $G_1$, $G_2$ in the family $\mathcal{G}$, there exists $G_3\in \mathcal{G}$ such that $G_3\subseteq \uparrow G_2\cap \uparrow G_1$.

The following two lemmas are well-known and can be easily verified.

\begin{lemma}\label{subspace-irr-1}
Let $X$ be a space and $Y$ a subspace of $X$. Then the following conditions are equivalent for a
subset $A\subseteq Y$:
\begin{enumerate}[\rm (1)]
	\item $A$ is an irreducible subset of $Y$.
	\item $A$ is an irreducible subset of $X$.
	\item ${\rm cl}_X A$ is an irreducible closed subset of $X$.
\end{enumerate}
\end{lemma}

\begin{lemma}\label{continuous-mapping-irr}
If $f: X\rightarrow Y$ is continuous and $A\in\ir(X)$, then $f(A)\in\ir(Y)$.
\end{lemma}

A subset $W$ of a $T_0$-space $X$ is called \emph{saturated} if $W$ equals the intersection of all open sets containing it (equivalently, $W\in \mathbf{up}(X)$). We shall use $Q(X)$ to denote the set of all nonempty compact saturated subsets of $X$ and endow it with the \emph{Smyth order} $\sqsubseteq$, that is, for $K_1,K_2\in Q(X)$, $K_1\sqsubseteq K_2$ if{}f $K_2\subseteq K_1$. The upper Vietoris topology on $Q(X)$ is the topology that has $\{\Box U: U\in \mathcal{O}(X)\}$ as a base, where $\Box U=\{Q\in Q(X):Q\subseteq U\}$, and the resulting space is called the \emph{Smyth power space} or \emph{upper space} of $X$ and is denoted by $P_S(X)$. For $A\in \Gamma(X)$, define $\Diamond A=\{Q\in Q(X):Q\cap A\neq\emptyset\}$. Clearly, $Q(X)\setminus \Box U=Q(X)\setminus \Diamond (X\setminus U$ for any $U\in\mathcal O(X)$.

\begin{lemma}\label{xi-continuous}
	For a $T_0$ space $X$, then the canonical mapping $\xi_X : X\longrightarrow P_S(X), x\mapsto\uparrow x$, is a dense topological embedding.
\end{lemma}
\begin{proof}
For $U\in\mathcal O(X)$, we have
$$\xi_X^{-1}(\Box U)=\{x\in X : ~\!\!\!\uparrow x\in\Box U\}=\{x\in X : ~\!\!\!\uparrow x\subseteq U\}=U, \hbox{~and}$$.
$$\hbox{~if~} U\neq\emptyset, \hbox{~then~} \emptyset\neq\{\uparrow x : x\in U\}=\xi_X (X)\cap\Box U.$$
\noindent So $\xi_X$ is continuous and $\xi_X (X)$ is dense in $P_S(X)$.
In addition, we have that $\xi_X(U)=\{\uparrow x : x\in U\}
=\{\uparrow x : ~\!\!\!\uparrow x\in\Box U\}
=\xi_X(X)\cap \Box U$,
\noindent which implies that $\xi_X$ is an open mapping to $\xi_X(X)$, as a subspace of $P_S(X))$.
As $X$ is $T_0$, $\xi_X$ is an injection. Thus $\xi_X$ is a dense topological embedding.
\end{proof}

A \emph{net} $(x_i)_{i\in I}$ in a set $X$ is a mapping from a directed set $I$ to $X$.  Suppose that $H(x)$ is a property of the elements $x\in X$. We say that $H(x)$ \emph{holds eventually} in the net $(x_i)_{i\in I}$ if there is a $i_0\in J$ such that $H(x_i)$ is true whenever $i_0\leq i$.

\section{Topological Rudin lemma}  %ÒÔÏÂÖ÷Òª½á¹û

Rudin's Lemma is a useful tool in topology and plays a crucial role in domain theory (see \cite{GHKLMS}). Rudin proved her lemma by transfinite methods. In \cite{HK}, Heckmann and Keimel presented the following topological variant of Rudin's Lemma.

\begin{lemma}\emph{(\cite[Rudin's Lemma]{GHKLMS})}
Let $\mathcal{F}$ be a directed family of nonempty finite subsets of a partially ordered set $P$. Then there exists a directed set $D\subseteq \bigcup_{F\in\mathcal{F}} F$ such that $D\cap F\neq \emptyset$ for all $F\in\mathcal{F}$.
\end{lemma}

\begin{lemma}\emph{(\cite[Lemma 3.1]{HK})}\label{T Rudin}  \emph{(Topological Rudin Lemma)} Let $X$ be a $T_0$-space and $\mathcal{A}\in \ir(P_S(X))$. Then every closed set $C {\subseteq} X$  that
meets all members of $\mathcal{A}$ contains a minimal irreducible closed subset $A$ that still meets all
members of $\mathcal{A}$.
\end{lemma}

Applying Lemma \ref{T Rudin} to the Alexandroff topology on a poset $P$, one obtains the original Rudin's Lemma (see \cite{Rudin}).

\begin{lemma}\label{rudin} \emph{(Rudin's Lemma)} Let $P$ be a poset, $C$ a nonempty lower subset of $P$ and $\mathcal F\subseteq \mathbf{Fin} P$ a filtered family with $C\cap \ua F\neq\emptyset$ for all $\ua F\in \mathcal F$. Then there exists a directed subset $D$ of $C$ such that $\ua F\cap \downarrow D\neq\emptyset$ for all $\ua F\in \mathcal F$.
\end{lemma}

On some occasions we only need the following consequence of Rudin's Lemma (cf. \cite[Lemma 2.4]{HK}).

\begin{corollary}\label{cor-rudin}  Let $P$ be a dcpo, $U$ a nonempty Scott open subset of $P$ and $\mathcal F\subseteq \mathbf{Fin} P$ a filtered family such that $\bigcap \mathcal F\subseteq U$ holds. Then $\ua F\subseteq U$ for some $\ua F\in \mathcal F$.
\end{corollary}

\begin{lemma}\emph{(\cite{WHX})}\label{Rudin irr} Let $X$ be a $T_0$-space and $\mathcal{F}\subseteq \mathbf{Fin}~\!\! X$. If $\mathcal F\in \ir(P_S(X))$,
 then there exists $A\in \ir(X) $ with $A\subseteq \bigcup \mathcal{F}  $ such that  $A\cap \ua F\neq \emptyset$ for any $\ua F\in \mathcal{F} $.
\end{lemma}

%ÒÔÉÏÊÇrudinÒýÀíÔÚ¼ÈÔ¼¼¯ÉϵÄÍƹã

\begin{lemma}\emph{(\cite{WHX})}\label{Rudin irr closed set} Let $X$ be a $T_0$-space $X$ and $\mathcal{K}\in \ir(P_S(X))$. If $B$ is a closed set that intersects every element of $\mathcal{K}$. Then $\{\ua (K\cap B) : K\in \mathcal{K}\}\in \ir (P_S(X))$.
\end{lemma}

Clearly, Lemma \ref{Rudin irr closed set} can be restated as follows.

\begin{lemma}\emph{(\cite{WHX})}\label{Rudin irr open set}  Let $X$ be a $T_0$-space $X$ and $\mathcal{K}\in \ir(P_S(X))$. If $U$  is an open set satisfying $K\setminus U\neq \emptyset$ for any  $K\in \mathcal{K}$.
 Then $\{\ua (K\setminus U) : K\in \mathcal{K}\}\in\ir (P_S(X))$.
\end{lemma}

For any $F\subseteq X$ and any upper set $U$ of $X$, it is straightforward to verify that $\ua (\ua F\setminus U)=\ua (F\setminus U)$. So by Lemma \ref{Rudin irr open set} we get the following.

\begin{corollary}\emph{(\cite{WHX})}\label{cor-Rudin irr open set}  Let $X$ be a $T_0$-space $X$, $U\in \mathcal O(X)$ and $\mathcal{F}\subseteq \mathbf{Fin}~\!\! X$. If $\mathcal F\in \ir(P_S(X))$ and $\ua F\setminus U\neq \emptyset$ for any  $\ua F\in \mathcal{F}$,
 then $\{\ua (F\setminus U) : \ua F\in \mathcal{F}\}=\{\ua (\ua F\setminus U) : \ua F\in \mathcal{F}\}\in\ir (P_S(X))$.
\end{corollary}

\begin{definition} (\cite[Definition 3.1]{ZH})\label{SI-open} %%%D. Zhao and W. K. Ho, ¡°On topologies defined by irreducible sets¡±, J. Log. Algebr. Methods Program. 84:1 (2015), 185¨C195.%%%%%
Let $X$ be a $T_0$-space. A subset $U$ of $X$ is called \emph{SI}-\emph{open} if the following two conditions are satisfied:
\begin{enumerate}[\rm (1)]
\item $U$ is an open set in $X$, and

\item for any $F\in \ir(X)$, $\bigvee F\in U$ implies $F\cap U\neq\emptyset$.
\end{enumerate}
The set of all SI-open sets in $X$ is denoted by $\mathcal O_{\SI}(X)$. It is straightforward to verify that $\mathcal O_{\SI}(X)$ is a topology on $X$, called the \emph{SI}-\emph{topology}. The space $(X, \mathcal O_{\SI}(X))$ will also be simply written as $SI(X)$. It is readily checked that for any $x\in X$  it holds that $\operatorname{cl}_X(\{x\})=\operatorname{cl}_{SI(X)}(\{x\})$.
\end{definition}

The following result can be considered as a generalization of Corollary \ref{cor-rudin}.

\begin{proposition}\emph{(\cite{WHX})}\label{rudin yunyong}
Let $X$ be an irreducible-complete $T_0$-space, $U$ be an SI-open set and $\mathcal{F}\subseteq \mathbf{Fin}~\!\! X$. If $\mathcal F\in \ir(P_S(X))$ and $\mathop{\bigcap}\mathcal{F} \subseteq U$, then  $\ua F\subseteq U$ for some $\ua F\in\mathcal{F}$.
\end{proposition}

\section{GIS-convergence}

\begin{definition}(\cite{CK}) \label{GELB}%%%%%%³ÂÇïÑ࣬¿Ü »Ô A characterization of quasi - continuous domains by nets    Journal of Sichuan University (Natural Science Edition)  doi£º10.3969/j.issn.0490-6756.2014.03.03  51¾í3ÆÚ
Let $P$ be a poset, $F\in P^{(<\omega)}$ and $(x_j)_{j\in J}$ be a net in $P$. $F$ is called a \emph{generalized eventual lower bound} of $(x_j)_{j\in J}$ if $x_j\in \ua F$ eventually, that is, there exists $k\in J$ such that $x_j\in \ua F$ for all $j\geq k$. Let $Elb(x_j)_{j\in J}$ denote the set of all generalized eventual lower bounds of the net $(x_j)_{j\in J}$.
\end{definition}

\begin{definition}(\cite[Definition 3.1]{ZL})\label{def:3.1}
Let $L$ be a dcpo, $x\in L$ and $(x_i)_{i\in I}$ a net on $L$.
$(x_i)_{i\in I}$ is said $GS$-convergent to $x$ iff there is a filtered family $\mathcal{F}=\{F:F \text{ is finite}\}$ such that for all $F\in\mathcal{F}$, $x_i\in\uparrow F$ eventually and $\bigcap_{F\in\mathcal{F}}\uparrow F\subseteq\uparrow x$.
It is denoted by $(x_i)_{i\in I}\stackrel{GS}\longrightarrow x$.
\end{definition}

\begin{definition}(\cite{ZLL})\label{IS-converges}
Let $(X,\tau)$ be a $T_0$ space, $x\in X$, and $(x_j)_{j\in J}$ be a net in $X$.
\begin{enumerate}[\rm (1)]
\item We say that $x$ is an eventual lower bound of $(x_j)_{j\in J}$ if $(x_j)_{j\in J}$ is eventually $\uparrow x$.

\item $(x_j)_{j\in J}$ is said to $IS$-converge to $x$, if there exists an irreducible subset $F$ of eventual lower bounds of $(x_j)_{j\in J}$ with $F^u\neq\emptyset$ such that $x\in \operatorname{cl}_{\mathcal{O}_{SI}}(F)$.
In this case, we write $(x_j)_{j\in J}\xrightarrow{IS} x$.
For a $T_0$ space $(X,\tau)$, it is obvious that
$
\mathcal O(\mathcal{IS}(X))=\bigl\{U\subseteq X\mid \text{whenever } (x_j)_{j\in J}\xrightarrow{IS} x \text{ and } x\in U,\text{ then }(x_j)_{j\in J}\text{ is eventually in }U\bigr\}
$
is a topology, called the $IS$-topology on $X$. The set $U$ in $\mathcal O(\mathcal{IS}(X))$ is called $IS$-open.
\end{enumerate}
\end{definition}

\begin{definition} \label{GIS-converges}
A net $(x_i)_{i\in I}$ is called \emph{GIS}-\emph{convergent} to a point $x$ in a~$T_0$-space $(X,\tau)$ if there exists an $\mathcal{F}\subseteq X^{(<\omega)}$ such that $\{\ua F : F\in\mathcal F\}\in \ir (P_S(X))$ and the following two conditions are satisfied:
\begin{enumerate}[\rm (i)]
\item $\mathcal{F}\subseteq Elb(x_j)_{j\in J}$;

\item for any $\mathcal{G}\in \mathcal{F}^{(<\omega)}$, $\bigcap\{\ua G: G\in \mathcal{G}\}\neq\emptyset$;

\item  $\ua x\in \cl_{P_S(X,\mathcal{O}_{SI}(X))}\{\ua F: F\in \mathcal{F}\}$.

\end{enumerate}
\noindent And in this case, we write $(x_i)_{i\in I}\stackrel{GIS}\longrightarrow x$. Let $\mathcal{GIS}(X)=\{((x_i)_{i\in I}, x) : (x_i)_{i\in I} \mbox{~is~a~net~in~} X, x\in X $ \mbox{~and~} $ (x_i)_{i\in I}\stackrel{GIS}\longrightarrow x\}$.
$\mathcal O(\mathcal{GIS}(X))=\bigl\{ U\,\big|\, \text{if } ((x_{i})_{i\in I},x)\in \mathcal{GIS}(X) \text{ and } x\in U,\text{ then the net }(x_{i})_{i\in I}\text{ is eventually in }U \bigr\}.$ The set $U$ in $\mathcal O(\mathcal{GIS}(X))$ is called $GIS$-open.
\end{definition}

\begin{remark}\label{GISremark} For a $T_0$-space $X$, we have the following conclusions:
\begin{enumerate}[\rm (1)]
\item The constant net $(x)_{i\in I}$ in $X$ with value $x$ GIS-converges to $x$. If $D$ is a directed set in $X$ and $x \in \operatorname{cl}_{\mathcal O_{SI}(X)}D$, then obviously the net consisting of the directed set $D$ itself GIS-converges to $x$.

\item If $(x_i)_{i\in I}\stackrel{GIS}\longrightarrow x$ in $X$, then $(x_i)_{i\in I}\stackrel{GIS}\longrightarrow y$ for any $y\leq x$. Thus the GIS-convergence points of a net are generally not unique.

\item Let $P$ be a dcpo. Then the GIS-convergence in $(P, \alpha(P))$ coincides with the
GS-convergence in $P$.

\item If $U\in \mathcal O(\mathcal{GIS}(X))$, then $U=\ua U$.

Suppose that $x\in U\in \mathcal O(\mathcal{GIS}(X))$ and $x\leq y$. Consider the constant net $(y)_{i\in I}$ in $X$ with value $y$. Clearly,
$(y)_{i\in I}\xrightarrow{GIS} y$, and hence $(y)_{i\in I}\xrightarrow{GIS} x$. As $x\in U\in\mathcal O(\mathcal{GIS}(X))$, $(y)_{i\in I}$ is eventually in $U$, that is, $y\in U$. So $U$ is an upper set.

\item If $(x_i)_{i\in I}\stackrel{IS}\longrightarrow x$ in $X$, then $(x_i)_{i\in I}\stackrel{GIS}\longrightarrow x$ in $X$; therefore, $\mathcal O(\mathcal{GIS}(X))\subseteq \mathcal O(\mathcal{IS}(X))$.

\item $\mathcal O_{SI}(X)\subseteq \mathcal O(\mathcal{GIS}(X))$ by the Definition \ref{IS-converges} and Definition \ref{GIS-converges}; thus, $\mathcal O(\mathcal{GIS}(X))$ is a $T_0$ spaces.

\end{enumerate}
\end{remark}

\begin{theorem}\label{GIS-eq}
Let $(X,\tau)$ be a $T_0$-space, $U\subseteq X$. Then the following statements are equivalent:
\begin{enumerate}[\rm (1)]
    \item $U\in \mathcal O(\mathcal{GIS}(X))$;
    \item  $U$ is an upper set. Let $x\in X$ and $\mathcal{F}\subseteq \mathbf{Fin}X$ satisfying $\mathcal{F}\in \ir (P_S(X))$. If for any $\mathcal{G}\in \mathcal{F}^{(<\omega)}$, $\bigcap\{\ua G: \ua G\in \mathcal{G}\}\neq\emptyset$ and $\ua x\in \cl_{P_S(X,\mathcal{O}_{SI}(X))}\{\ua F: \ua F\in \mathcal{F}\}\cap \Box U\neq\emptyset$, then $\bigcap\{\ua G: \ua G\in \mathcal{G}_0\}\subseteq U$ for some $\mathcal{G}_0\in \mathcal{F}^{(<\omega)}$.

\end{enumerate}
\end{theorem}

\begin{proof}
(1) $\Rightarrow$ (2):
 By Remark 4.5(4),  we have $U=\ua U$.
Let $U\in \mathcal O(\mathcal{GIS}(X))$,  $x\in X$ and let $\mathcal{F}\subseteq \mathbf{Fin}X$ satisfying $\mathcal{F}\in \ir (P_S(X))$ such that  for any $\mathcal{G}\in \mathcal{F}^{(<\omega)}$ we have $\bigcap\{\ua G: \ua G\in \mathcal{G}\}\neq\emptyset$ and $\ua x\in \cl_{P_S(X,\mathcal{O}_{SI}(X))}\{\ua F: \ua F\in \mathcal{F}\}\cap \Box U\neq\emptyset$. Let $D_{\mathcal{F}}=\{\bigcap\mathcal{G}\mid \mathcal{G} \text{ is a nonempty finite subset of }\mathcal{F}\}$ and $D_{\mathcal{F}}^X=\{(d,D)\in X\times D_{\mathcal{F}}\mid d\in D\}$.
 The order $\leq$ on $D_{\mathcal{F}}^X$ is defined as follows: $\forall\,(d_1,D_1),(d_2,D_2)\in D_{\mathcal{F}}^X,\;(d_1,D_1)\leq (d_2,D_2)\iff D_2\subseteq D_1.$ Then $(D_{\mathcal{F}}^X,\leq)$ is a preordered set, and thus $D_{\mathcal{F}}^X$ is a directed set.
Let $x_{(d,D)}=d$ for any $(d,D)\in D_{\mathcal{F}}^X$. Then the net $(x_{(d,D)})_{(d,D)\in D_{\mathcal{F}}^X}\xrightarrow{GIS} x$. It follows from $U\in \mathcal O(\mathcal{GIS}(X))$ that $(x_{(d,D)})_{(d,D)\in D_{\mathcal{F}}^X}$ is eventually in $U$, that is, there exists $(d_0,D_0)\in D_{\mathcal{F}}^X$ such that $x_{(d,D)}\in U$ for any $(d,D)\geq (d_0,D_0)$. Specially, $D_0\subseteq U$. Hence, there exists a nonempty finite subset $\mathcal{G}_0$ of $\mathcal{F}$ such that $\bigcap\{\ua G: \ua G\in \mathcal{G}_0\}\subseteq U$.

(2) $\Rightarrow$ (1): Let $(x_j)_{j\in J}\xrightarrow{GIS} x$ and $x\in U$. Then there exists an irreducible set $\mathcal{F}\subseteq \mathbf{Fin}X$ such that for any $\mathcal{G}\in \mathcal{F}^{(<\omega)}$ we have $\bigcap\{\ua G: \ua G\in \mathcal{G}\}\neq\emptyset$ and $\ua x\in \cl_{P_S(X,\mathcal{O}_{SI}(X))}\{\ua F: \ua F\in \mathcal{F}\}$, and for any $\ua F\in \mathcal{F}$, $x_j\in\ua F$ eventually. Then $\ua x\in \cl_{P_S(X,\mathcal{O}_{SI}(X))}\{\ua F: \ua F\in \mathcal{F}\}\cap \Box U$. It follows from the hypothesis that there exists a nonempty finite subset $\mathcal{G}_0$ of $\mathcal{F}$ such that $\bigcap\{\ua G: \ua G\in \mathcal{G}_0\}\subseteq U$. Since $\mathcal{G}_0$ is a nonempty finite set, there exists $i_0\in J$ such that $x_j\in \bigcap\{\ua G: \ua G\in \mathcal{G}_0\}$ for any $j\geq i_0$.
So we conclude that $(x_j)_{j\in J}$ is eventually in $U$. Therefore, $U\in\mathcal O(\mathcal{GIS}(X))$.
\end{proof}

\begin{proposition}\label{GIS-eq-1}
Let $(X,\tau)$ be a $T_0$-space, $U\subseteq X$. If $U$ satisfies the following conditions:
\begin{enumerate}[\rm (1)]

\item $U$ is an upper set;

\item Let $x\in X$ and $\mathcal{F}\subseteq \mathbf{Fin}X$ satisfying $\mathcal{F}\in \ir (P_S(X))$. If for any $\mathcal{G}\in \mathcal{F}^{(<\omega)}$, $\bigcap\{\ua G: \ua G\in \mathcal{G}\}\neq\emptyset$ and $\ua x\in \cl_{P_S(X,\mathcal{O}_{SI}(X))}\{\ua F: \ua F\in \mathcal{F}\}\cap \Box U$, then $\mathcal{F}\cap \Box U\neq\emptyset$;
\end{enumerate}
 then $U$ is a $GIS$-open set.
\end{proposition}

\begin{proof}
Let $x\in X$ and $\mathcal{F}\subseteq \mathbf{Fin}X$ satisfying $\mathcal{F}\in \ir (P_S(X))$ such that for any $\mathcal{G}\in \mathcal{F}^{(<\omega)}$, $\bigcap\{\ua G: \ua G\in \mathcal{G}\}\neq\emptyset$ and $\ua x\in \cl_{P_S(X,\mathcal{O}_{SI}(X))}\{\ua F: \ua F\in \mathcal{F}\}\cap \Box U$. Then by (2) there exists $\ua G\in \mathcal{F}\cap \Box U$. Since $U$ is an upper set, we have $\ua G\subseteq U$. It follows from Theorem~\ref{GIS-eq} that $U$ is a $GIS$-open set.
\end{proof}

\begin{remark}
The converse of Proposition \ref{GIS-eq-1} may not be true. Please see the following Example.
\end{remark}

\begin{example}
Let $P=\{a,c\}\cup\{b_n\mid n\in\mathbb{N}^+\}$, where $\mathbb{N}^+$ is the set of all positive integers. The order $\leq$ on $P$ is defined as follows (see Fig.~2):

\begin{figure}[htbp]
\centering
\includegraphics[width=0.45\textwidth]{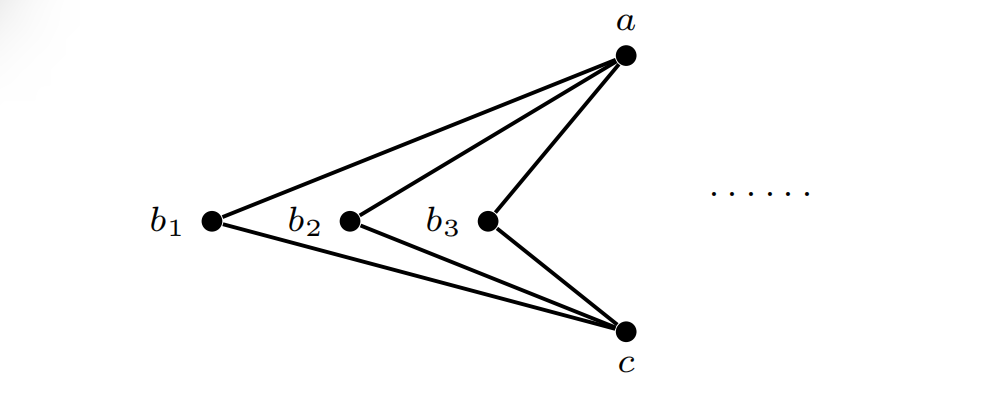}
\caption{The poset $P$ in Example~4.9}
\label{fig:poset_example4.9}
\end{figure}

\begin{itemize}
\item[(i)] $c\leq b_n\leq a$ for all $n\in\mathbb{N}^+$;
\item[(ii)] For any $m\neq n\in\mathbb{N}^+$, $b_m$ and $b_n$ are incomparable.
\end{itemize}
Considering the topology $\tau=\upsilon(P)$ on $P$. We claim that $\ua b_1\in \mathcal O(\mathcal{GIS}(P,\tau))$.
Let $\mathcal{F}\subseteq \mathbf{Fin}P$ satisfying $\mathcal{F}\in \ir (P_S(P,\tau))$ such that for any $\mathcal{G}\in \mathcal{F}^{(<\omega)}$, $\bigcap\{\ua G: \ua G\in \mathcal{G}\}\neq\emptyset$ and $\emptyset\neq\cl_{P_S(X,\mathcal{O}_{SI}(P))}\{\ua F: \ua F\in \mathcal{F}\}\cap \Box \ua b_1\subseteq \{\ua b_1, \ua a\}$. If $\ua a\in \mathcal{F}$, then $\ua a\in\{\ua F: \ua F\in \mathcal{F}\}\cap\Box\ua b_1$. If $\ua a\notin \mathcal{F}$, then $\mathcal{F}\subseteq \mathbf{Fin}(\{c\}\cup\{b_n\mid n\in\mathbb{N}^+\})$. Suppose that $\mathcal{F}=\{\ua c\}$, $\cl_{P_S(P,\mathcal{O}_{SI}(P))}\{\ua F: \ua F\in \mathcal{F}\}=\{\ua c\}$ which contradicts the fact that $\cl_{P_S(P,\mathcal{O}_{SI}(P))}\{\ua F: \ua F\in \mathcal{F}\}\cap \Box \ua b_1\neq\emptyset$. Then $\mathcal{F}\cap\{\ua G:G\in(\{b_n\mid n\in\mathbb{N}^+\})^{(<\omega)}\}\neq\emptyset$.
Suppose that  $\bigcap\mathcal{F}\nsubseteq \ua b_1$. There exists $\ua G_0\in \{\ua G:G\in(\{b_n\mid n\in\mathbb{N}^+\})^{(<\omega)}\}\setminus\{\ua b_1\}$ such that $\bigcap\mathcal{F}=\ua G_0$. Hence, $\cl_{P_S(P,\mathcal{O}_{SI}(P))}\{\ua F: \ua F\in \mathcal{F}\}\subseteq \downarrow_{Q(P,\upsilon(P))}\ua G_0$ and $\downarrow_{Q(P,\upsilon(P))}\ua G_0\cap \Box \ua b_1=\emptyset$ which contradicts the fact that $\cl_{P_S(P,\mathcal{O}_{SI}(P))}\{\ua F: \ua F\in \mathcal{F}\}\cap \Box \ua b_1\neq\emptyset$.
Thus $\bigcap\mathcal{F}\subseteq \ua b_1$, Hence, there exists $\mathcal{G}_0\in \mathcal{F}^{(<\omega)}$ such that $\bigcap\mathcal{G}_0\subseteq \ua b_1$.
By Theorem 4.6, we show that $\ua b_1\in \mathcal O(\mathcal{GIS}(P,\tau))$.

Let $F=\{b_m\mid m\in\mathbb{N}^+, m\geq 2\}$. Then $F$ is an irreducible set and $\{\ua G: G\in  F^{(<\omega)}\}\in \ir (P_S(P,\tau))$, and thus $\bigcap\{\ua G: G\in  F^{(<\omega)}\}=\ua a$ and $(P,\tau)$ is an irreducible-complete $T_0$-space. So $\ua a\in \cl_{P_S(P,\mathcal{O}_{SI}(P))}\{\ua G: G\in  F^{(<\omega)}\}$ by Proposition 3.10 and thus $\ua a\in \cl_{P_S(P,\mathcal{O}_{SI}(P))}\{\ua G: G\in  F^{(<\omega)}\}\cap \Box \ua b_1$. Obviously, $\{\ua G: G\in  F^{(<\omega)}\}\cap\Box\ua b_1=\emptyset$. So we conclude that the converse of Proposition~\ref{GIS-eq-1} may not be true.
\end{example}

\begin{remark}
\begin{enumerate}[\rm (1)]
\item A $GIS$-open set may not be an $SI$-open set. For example, let $P$ be the $T_0$ space in Example 4.9. Then $\ua b_1$ is a $GIS$-open set. Since $\ua b_1 \notin \upsilon(P)=\tau$, we have that $\ua b_1 \notin \mathcal{O}_{SI}$.

\item Neither $\tau$ nor $\mathcal O(\mathcal{GIS}(X))$ is contained in the other. For example, let $P$ be the $T_0$ space as in Example 4.9. By (1), we have that $\uparrow b_1\notin\tau$ and $\uparrow b_1\in\mathcal O(\mathcal{GIS}(P,\upsilon(P)))$. let $L=\{\top\}\cup \mathbb{N}^+$, where $\top$ is the greatest element of $L$ and $\mathbb{N}^+$ is the set of all positive integers. Endow $L$ with the Alexandroff topology $\alpha(L)$. Then $\{\top\}\in\alpha(L)$ and $SI(L)=(L,\sigma(L))$. let $\mathcal{F}=\{\ua n: n\in \mathbb{N}^+\}$. Actually,  $\bigcap\mathcal{F}=\ua \top$ and $\ua \top\in \cl_{P_S(L,\mathcal{O}_{SI}(L))}\mathcal{F}$. But $\mathcal{F}\cap\Box \ua \top=\emptyset$. Hence, $\{\top\}$ is not a $GIS$-open set.

    \end{enumerate}
\end{remark}

The following example show that a Scott open set may not be an $GIS$-open set.

\begin{example}

Let $L = \mathbb{N}^{+} \cup \{a,b,c\}$. The order $\leq$ on $L$ is defined as follows:
\begin{enumerate}[\rm (i)]
    \item $a \leq b$, $a \leq c$;
    \item $\forall\, n \in \mathbb{N}^{+} ,\ n \leq b$ and $n \leq c$;
    \item $\forall\, m,n \in \mathbb{N}^{+} ,\ m \leq n$ if and only if they satisfy the usual order on $\mathbb{N}^{+} $.
\end{enumerate}

\begin{figure}[htbp]
\centering
\includegraphics[width=0.25\textwidth]{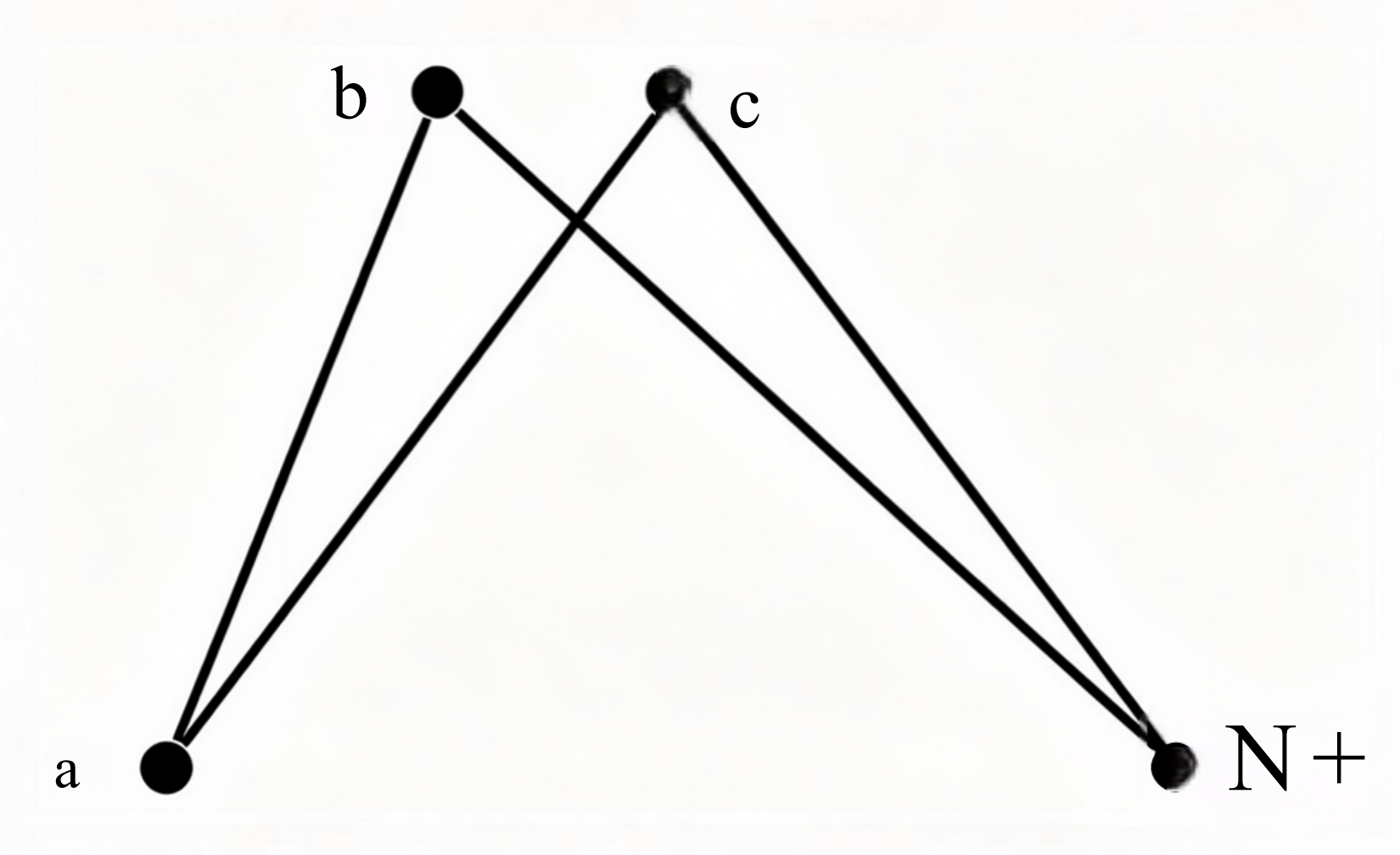}
\caption{The poset $L$ in Example~4.11}
\label{fig:poset_example4.9}
\end{figure}

Then $(L,\leq)$ is a poset (see Fig. 2). Considering the topology $ (L,\upsilon(L))$.  Considering the set $A=\{a,b,c\}$. Actually, $A$ is a Scott open set, $\mathcal{F}=\{\ua n: n\in \mathbb{N}^+\}\in \ir(P_S(L,\upsilon(L)))$, and $\bigcap\mathcal{F}=\{b,c\}$.  As $\cl_{\upsilon(L)}(\mathbb{N}^+) = \mathbb{N}^+ \cup \{a\}$, then $\ua a\in\cl_{P_S((L,\upsilon(L)))}\mathcal{F}\subseteq\cl_{P_S(L,\mathcal{O}_{SI}(L))}\mathcal{F}$. But $\mathcal{F}\cap\Box A=\emptyset$. Hence, $A$ is not  a $GIS$-open set.
\end{example}

\begin{definition}(\cite{E2009})
Let $(X,\tau)$ be a $T_0$-space.

\begin{enumerate}[\rm (1)]
\item The space $(X,\tau)$ is called $C$-space if for any $U\in\tau$ and any $x\in U$, there exists $y\in X$ such that $x\in \ii_\tau \ua y \subseteq \ua y \subseteq U$.
\item The space $(X,\tau)$ is called \textit{locally hypercompact} if for any $U\in\tau$ and any $x\in U$, there exists $F\in X^{(<\omega)}$ such that $x\in \ii_\tau \ua F \subseteq \ua F \subseteq U$.

\end{enumerate}

\end{definition}

\begin{corollary}
Let $(X,\tau)$ be a $T_0$ space. If $SI(X)$ is a  locally hypercompact space, then $\mathcal{O}_{SI}(X) = \mathcal O(\mathcal{GIS}(X))$.
\end{corollary}

\begin{proof}
By Remark 4.5(6), it suffices to prove that $\mathcal{O}_{SI}(X) \supseteq \mathcal O(\mathcal{GIS}(X))$. Let $x\in X$ and $U\in \mathcal O(\mathcal{GIS}(X))$ with $x\in U$. As $SI(X)$ is a  locally hypercompact space, $\mathcal{F}=\{\ua F:x\in \ii_{\mathcal{O}_{SI}(X)}\ua F~\text{and}~F\in X^{(<\omega)}\}$ is filtered. Let $I=\{(\ua F,e)\in \mathcal F\times X: e\in \ua F\}$ and define an order $\leq$ on $I$ as follows:  $(\ua F_1,e_1)\leq (\ua F_2,e_2)$ iff $\ua F_2\subseteq \ua F_1$.  For any $(\ua F,e)\in I$, $x_{(\ua F,e)}=e$. It is easy to see that $\mathcal{F}$ and $(x_i)_{i\in I}$ satisfy conditions (i)--(iii) of Definition 4.4. Therefore,  $(x_i)_{i\in I}\stackrel{GIS}\longrightarrow x$. From $x\in U\in \mathcal O(\mathcal{GIS}(X))$, it follows that $x_i\in U$ eventually, that is, there exists $\ua G\in\mathcal F$ such that $x\in \ii_{\mathcal{O}_{SI}(X)} G\subseteq \ua G\subseteq U$. Therefore, $\mathcal{O}_{SI}(X) \supseteq \mathcal O(\mathcal{GIS}(X))$.
\end{proof}

\begin{corollary}
Let $(X,\tau)$ be a $d$-space. If $SI(X)$ is a  locally hypercompact space, then $\mathcal{O}_{SI}(X) = \mathcal O(\mathcal{GIS}(X))=\tau$.
\end{corollary}

\begin{proof}
By Corollary 4.13, it suffices to show that $ \mathcal O(\mathcal{GIS}(X))\supseteq\tau$. Let $U\in \tau$, $x\in U$ and $(x_i)_{i\in I}$ with $(x_i)_{i\in I}\stackrel{GIS}\longrightarrow x$. Then there exists $\mathcal{F} \subseteq X^{(<\omega)}$ such that $\{\uparrow G : G \in \mathcal{F}\} \in \ir(P_S(X))$ and conditions (i)--(iii) in Definition~4.4 are satisfied. Since $SI(X)$ is a  locally hypercompact space, $\mathcal{G}=\{\ua G: x\in \ii_{\mathcal{O}_{SI}(X)}\ua G~\text{and}~G\in X^{(<\omega)}\}$ is filtered with $\mathcal{G}\subseteq \da_{Q(X)} \{\ua G : G \in \mathcal{F}\} $ and $\ua x\in \cl_{P_S(SI(X))}\mathcal{G}$.  Suppose, for contradiction, that $\ua G\nsubseteq U$ for any $\ua G\in\mathcal{G}$. Then $\{\ua(\ua G\backslash U): \ua G\in\mathcal{G}\}$ is filtered and $\ua(\ua G\backslash U)=\ua(G\backslash U)$ for any $\ua G\in\mathcal{G}$. By Lemma 3.1, there exists a directed set  $D\subseteq \bigcup_{\ua G\in \mathcal{G}}G\backslash U$ such that $D\cap (G\backslash U)\neq\emptyset$. As $SI(X)$ is a  locally hypercompact space, it follow that $x\in \cl_{\mathcal{O}_{SI}(X)}D$. Since $(X,\tau)$ is a $d$-space, $x\in \cl_{\mathcal{O}_{SI}(X)}D=\cl_{\tau}D$ and hence $U\cap D\neq\emptyset$, which is a contradiction to $D\subseteq \bigcup_{\uparrow G\in \mathcal{G}}G\backslash U$. So there is $\ua G\in \mathcal{G}$ such that $x\in \ii_{\mathcal{O}_{SI}(X)} \ua G \subseteq \ua G \subseteq U$. Therefore, $ \mathcal O(\mathcal{GIS}(X))\supseteq\tau$.
\end{proof}

\begin{proposition} Let $(X, \tau)$ be a $T_0$ space. Then the following statements hold:

\begin{enumerate}[\rm (1)]
  \item $(X, \mathcal O(\mathcal{GIS}(X)))$ is a $T_0$ space;
  \item For any directed subset $D$ of $X$, $\operatorname{cl}_{\mathcal{O}_{SI}(X)}(D) = \operatorname{cl}_{\mathcal O(\mathcal{GIS}(X))}(D)$;
  \item The specialization orders of spaces $(X, \tau)$ and $(X, \mathcal O(\mathcal{GIS}(X)))$ coincide.
\end{enumerate}

\end{proposition}

\begin{proof}
(1) As for any $y\in X$ it follow that $X\backslash \da y=X\backslash \cl_{\mathcal{O}_{SI}(X)} y\in \mathcal{O}_{SI}(X)\subseteq \mathcal O(\mathcal{GIS}(X))$, $(X, \mathcal O(\mathcal{GIS}(X)))$ is a $T_0$ space;

(2)  Let $D$ be a directed subset of $X$. By $\mathcal{O}_{SI}(X)\subseteq \mathcal O(\mathcal{GIS}(X))$, it follow that $\cl_{\mathcal O(\mathcal{GIS}(X))}D\subseteq \cl_{\mathcal{O}_{SI}(X)}D$. Let $x\in \cl_{\mathcal{O}_{SI}(X)}D$. It is easy to know that $D\stackrel{GIS}\longrightarrow x$ and thus $x\in \cl_{\mathcal O(\mathcal{GIS}(X))}D$.

(3)  By (2), it follow that $\cl_{\mathcal O(\mathcal{GIS}(X))}\{x\}=\cl_{\mathcal{O}_{SI}(X)}\{x\}=\cl_{\tau}\{x\}$ and thus the specialization orders of spaces $(X, \tau)$ and $(X, \mathcal O(\mathcal{GIS}(X)))$ coincide.
\end{proof}

\begin{proposition}

Let $(X, \tau)$ be a $T_0$-space. Then the following statements are equivalent:

\begin{enumerate}[\rm (1)]
  \item $SI(X)$ is a locally hypercompact space;
  \item For any net $(x_i)_{i\in I}$ in $X$ and $x \in X$, $(x_i)_{i\in I} \xrightarrow{GIS} x$ if and only if $(x_i)_{i\in I}$ converges to $x$ with respect to the topology $\mathcal{O}_{SI}(X)$.
\end{enumerate}

\end{proposition}

\begin{proof}
(1) $\Rightarrow$ (2):  Suppose that $(x_i)_{i\in I} \xrightarrow{GIS} x$. There exists $\mathcal{F}\subseteq X^{(<\omega)} $ such that $\{\uparrow G : G \in \mathcal{F}\} \in \ir(P_S(X))$ and conditions (i)--(iii) in Definition~4.4 are satisfied. Let $U\in \mathcal{O}_{SI}(X)$ and $x\in X$ with $x\in U$. Thus, There exists $G_0\in \mathcal{F}$ such that $\ua G_0\subseteq U$.
As $G_0$ is a generalized eventual lower bound of $(x_i)_{i\in I}$, it follow that $x_i\in U$ eventually, and thus $(x_i)_{i\in I}$ converges to $x$ with respect to the topology $\mathcal{O}_{SI}(X)$.

Conversely, suppose that $(x_i)_{i\in I}$ converges to $x$ with respect to the topology $\mathcal{O}_{SI}(X)$. Since $SI(X)$ is a locally hypercompact space, then $\{\ua F:x\in \ii_{\mathcal{O}_{SI}(X)}\ua F~\text{and}~F\in X^{(<\omega)}\}$ is filtered. Set $\mathcal{G}=\{\ua F:x\in \ii_{\mathcal{O}_{SI}(X)}\ua F~\text{and}~F\in X^{(<\omega)}\}$. It is readily verified that $\mathcal{G}\subseteq Elb(x_i)_{i\in I}$ and $\ua x\in \cl_{\mathcal{O}_{SI}(X)}\mathcal{G}$, and thus $(x_i)_{i\in I} \xrightarrow{GIS} x$.

(2) $\Rightarrow$ (1): Set $\mathcal{N}(x)=\{U\in \mathcal{O}_{SI}(X): x\in U\}$. Let $I=\{(U,e)\in \mathcal{N}(x)\times X: e\in U\}$ and define an order $\leq$ on $I$ as follows:
$(U_1,e_1)\leq (U_2,e_2)$ iff $U_2\subseteq U_1$. Clearly, $I$ is a directed set. For each $(U,e)\in I$, let $x_{(U,e)} = e$.
Clearly, $(x_{(U,e)})_{(U,e)\in I}$ converges to $x$ with respect to the topology $\mathcal{O}_{SI}(X)$. By condition (2), we have $(x_i)_{i\in I} \xrightarrow{GIS} x$. Then exists $\mathcal{K}\subseteq X^{(<\omega)} $ such that $\{\ua K : K \in \mathcal{K}\} \in \ir(P_S(X))$ and conditions (i)--(iii) in Definition~4.4 are satisfied. Thus, $\ua x\in \cl_{P_S(X,\mathcal{O}_{SI}(X))}\{\ua K : K \in \mathcal{K}\}$. Suppose that $V\in \mathcal{N}(x)$. Thus, there is $K_0\in \mathcal{K}$ such that $\ua K_0\subseteq V$. As $K_0$ is a generalized eventual lower bound of $(x_i)_{i\in I}$, there is $(V_0,e_0)\in I$ such that $x_{(V_1,e_1)}=e_1\in \ua K_0$ for any $(V_0,e_0)\leq(V_1,e_1)$. Hence, $x\in V_0\subseteq \ua K_0\subseteq V$. Therefore, $SI(X)$ is a locally hypercompact space.
\end{proof}

\begin{definition}

A $T_0$-space $(X, \tau)$ is called a $GID_C$-space if
for any irreducible subfamily $\mathcal{F}$ of $P_S(X)$ with $\bigcap\mathcal{F} \neq \emptyset$ and $\mathcal{F}\subseteq \mathbf{Fin}~\!\!X$, there exists a filtered subfamily $\mathcal{G}$ with $\mathcal{G}\subseteq \mathbf{Fin}~\!\!X\cap \da_{Q(X)}\mathcal{F}$ such that $\operatorname{cl}_{P_S(X,\mathcal{O}_{SI}(X))}(\mathcal{G}) = \operatorname{cl}_{P_S(X,\mathcal{O}_{SI}(X))}(\mathcal{F})$.

\end{definition}

\begin{proposition}
Let $(X, \tau)$ be an $GID_C$-space, $x\in X$ and $\mathcal{F}\in \ir (P_S(X))$ with $\mathcal{F}\subseteq \mathbf{Fin}~\!\!X$ and $\bigcap\mathcal{F}\neq\emptyset$. If $\ua x\in \cl_{P_S(X,\mathcal{O}_{SI}(X))}\mathcal{F}$, then  there exists a net $(x_j)_{j\in J}$ in $\bigcup\mathcal{F}$ such that $(x_j)_{j\in J} \xrightarrow{GIS} x$.
\end{proposition}

\begin{proof}
 Let $\mathcal{F}$ be an irreducible set in $P_S(X)$ such that $\mathcal{F}\subseteq \mathbf{Fin}~\!\!X$, $\bigcap\mathcal{F}\neq\emptyset$ and $\ua x\in \cl_{P_S(X,\mathcal{O}_{SI}(X))}\mathcal{F}$. Since $X$ is an $GID_C$-space, there exists filtered subfamily $\mathcal{G}$ with $\mathcal{G}\subseteq \mathbf{Fin}~\!\!X\cap \da_{Q(X)}\mathcal{F}$ such that $\operatorname{cl}_{P_S(X,\mathcal{O}_{SI}(X))}(\mathcal{G}) = \operatorname{cl}_{P_S(X,\mathcal{O}_{SI}(X))}(\mathcal{F})$. Then $\ua x \in \cl_{\mathcal{O}_{SI}(X)}(\mathcal{G})$. It follows from $\bigcap\mathcal{F} \neq \emptyset$ that $\bigcap\mathcal{G} \neq \emptyset$.
Let $I=\{(\ua G,e)\in \mathcal{G}\times X: e\in \ua G\}$ and define an order $\leq$ on $I$ as follows:
$(\ua G_1,e_1)\leq (\ua G_2,e_2)$ iff $\ua G_2\subseteq \ua G_1$.  Clearly, $I$ is a directed set. For each $(\ua G,e)\in I$, let $x_{(\ua G,e)} = e$. Then $(x_i)_{i\in I}$ is a net in $\bigcup\mathcal{F}$, and thus $(x_i)_{i\in I} \xrightarrow{GIS} x$.
\end{proof}

\begin{theorem}
Let \((X, \tau)\) be an \(GID_C\)-space, \(U \subseteq X\). Then the following statements are equivalent:
\begin{enumerate}[\rm (1)]
    \item $U\in \mathcal{O}(\mathcal{GIS}(X))$;
   \item The set $U$ satisfies the following conditions:
    \begin{enumerate}[\rm (i)]
        \item $U=\ua U$;
        \item For any $\mathcal{F}\in \mathbf{Fin}~\!\!X$ satisfying $\mathcal{F}\in \ir(P_s(X))$ . If $\cl_{P_S(X,\mathcal{O}_{SI}(X))}\mathcal{F}\cap \square U\neq \emptyset$, then $\mathcal{F}\cap \square U\neq \emptyset$.
    \end{enumerate}
    \item The set $U$ satisfies the following conditions:
    \begin{enumerate}[\rm (i)]
        \item $U=\ua U$;
        \item For any directed set $D$, if $\cl_{\mathcal{O}_{SI}(X)}D\cap U\neq \emptyset$, then $U\cap D\neq \emptyset$.
    \end{enumerate}
\end{enumerate}
\end{theorem}

\begin{proof}

(1) $\Rightarrow$ (2):
Obviously, condition (1) implies that $U$ is an upper set. Moreover, since $X$ is a $GID_C$-space, there exists a filtered family $\mathcal{G}$ with $\mathcal{G}\subseteq \mathbf{Fin}~\!\!X\cap \da_{Q(X)}\mathcal{F}$ such that $\operatorname{cl}_{P_S(X,\mathcal{O}_{SI}(X))}(\mathcal{G}) = \operatorname{cl}_{P_S(X,\mathcal{O}_{SI}(X))}(\mathcal{F})$.
Suppose, for the sake of contradiction, that $\mathcal{G}\cap \square U= \emptyset$. Let $\mathcal{K}=\{\ua (\ua G\setminus U):\ua G\in \mathcal{G}\}$. It is obvious that $\mathcal{K}$ is a filtered family.
Furthermore, by Rudin's lemma, there exists a directed set $D$ with $D\subseteq\bigcup \{min(\ua K):\ua K\in \mathcal{K}\}$   such that   $D\cap min(\ua K)\neq\emptyset $ hold for any $\ua K\in \mathcal{K}$. Therefore, $\diamond\cl_{\mathcal{O}_{SI}(X)}D\supseteq \operatorname{cl}_{P_S(X,\mathcal{O}_{SI}(X))}(\mathcal{G})$. Thus, $\diamond\cl_{\mathcal{O}_{SI}(X)}D\cap\Box U\neq\emptyset$ by $\cl_{P_S(X,\mathcal{O}_{SI}(X))}\mathcal{F}\cap \square U\neq \emptyset$. Hence, $\cl_{\mathcal{O}_{SI}(X)}D\cap U\neq\emptyset$.  Pick an arbitrary element $x \in \cl_{\mathcal{O}_{SI}(X)}D\cap U$. By Remark 4.5(1), it follows that $D$ GIS-converges to $x$ and $x \in U$. By $U\in \mathcal{O}(\mathcal{GIS}(X))$, $D\cap U\neq\emptyset$ which is a contradiction. Therefore, $\mathcal{G}\cap \square U\neq \emptyset$. Consequently, $\mathcal{F}\cap \square U\neq \emptyset$.

(2) $\Rightarrow$ (1): By Definition 4.4, it is obvious that $U\in \mathcal{O}(\mathcal{GIS}(X))$.

(2) $\Rightarrow$ (3): It is obvious.

(3) $\Rightarrow$ (2):  Suppose $x \in U$ and $(x_i)_{i \in I}$ GIS-converges to $x$. By Definition 4.4, there exists an $\mathcal{F}\subseteq X^{(<\omega)}$ such that $\{\ua F : F\in\mathcal F\}\in \ir (P_S(X))$ and $\ua x\in \cl_{P_S(X,\mathcal{O}_{SI}(X))}\{\ua F: F\in \mathcal{F}\}$. In a similar way to the proof of $(1) \Rightarrow (2)$, we can prove that $U\in \mathcal{O}(\mathcal{GIS}(X))$.
\end{proof}

\section{GIS-continuous spaces}

In this section, we shall introduce the notion of $GIS$-continuous spaces and prove that a $T_0$ space is a $GIS$-continuous space if and only if the $GIS$-convergence is topological.

Let $(X, \tau)$ be a $T_0$ space. For any $F, G \in X^{(<\omega)}$, define $F \ll_{GIS} G$ if for any net $(x_j)_{j \in J}$ $GIS$-converges to $y\in G$, $ x_i\in \ua F$ holds eventually. We denote the set $\{F \in X^{(<\omega)} \mid F\ll_{GIS} x\}$ by $w_{GIS} (x)$, and the set $\{x \mid F \ll_{GIS} x\}$ by $\Uparrow_{GIS} F$.

\begin{remark}
Let $(X, \tau)$ be a $T_0$ space, $G, H, K, F\in X^{(<\omega)}$ and $x \in X$. Then the following statements hold:
\begin{enumerate}[(1)]
    \item If $F \ll_{GIS} x$, then $x \in \ua F$.

    \item If $F \ll_{GIS} G$, then $G \subseteq \ua F$.

    \item If $F \ll_{GIS} G$, then  $F \ll_{GIS} g$ for any $g\in G$.

    \item If $G \le H \ll_{GIS} K \le F$, then $G \ll_{GIS} F$.

    \item If $x \in \operatorname{int}_{\mathcal{O}(\mathcal{GIS}(X)}(\ua F)$, then $F \ll_{GIS} x$.

\end{enumerate}
\end{remark}

\begin{theorem}
Let $(X, \tau)$ be a $T_0$ space, $F\in X^{(<\omega)}$, $x \in X$. Consider the following statements:
\begin{enumerate}[(1)]
    \item $F \ll_{GIS} x$;

    \item If $\mathcal{F}\subseteq X^{(<\omega)}$ with $\{\ua F : F\in\mathcal F\}\in \ir (P_S(X))$ and the following two conditions are satisfied:

      \begin{enumerate}[\rm (i)]
        \item for any $\mathcal{G}\in \mathcal{F}^{(<\omega)}$, $\bigcap\{\ua G: G\in \mathcal{G}\}\neq\emptyset$;
         \item $\ua x\in \cl_{P_S(X,\mathcal{O}_{SI})}\{\ua F: F\in \mathcal{F}\}$.
       \end{enumerate}

        There exists an $\mathcal{G}_0\subseteq \mathcal{F}^{(<\omega)}$ such that $\bigcap\{\ua G: G\in \mathcal{G}_0\}\subseteq \ua F$.

    \item For any irreducible subset $F$ of $X$ with $F^{\ua} \neq \emptyset$, $x \in \operatorname{cl}_{\mathcal{O}_{SI(X)}}(F)$ implies there exists a nonempty finite subset $F_0$ of $F$ such that $F_0^{\ua} \subseteq \ua F$.
\end{enumerate}

Then we have $(1) \Leftrightarrow (2) \Rightarrow (3)$.

\end{theorem}

\begin{proof}

(1) $\Rightarrow$ (2): Let $F \ll_{GIS} x$, and let $\mathcal{F}\subseteq X^{(<\omega)}$ with $\{\ua F : F\in\mathcal F\}\in \ir (P_S(X))$ such that for any $\mathcal{G}\in \mathcal{F}^{(<\omega)}$, $\bigcap\{\ua G: G\in \mathcal{G}\}\neq\emptyset$ and $\ua x\in \cl_{P_S(X,\mathcal{O}_{SI}(X))}\{\ua F: F\in \mathcal{F}\}$.
It follows from the proof of Theorem 4.6 that there exists a net $(x_{(d,D)})_{(d,D) \in D_{\mathcal{F}}^X} \xrightarrow{GIS} x$.
Then $(x_{(d,D)})_{(d,D) \in D_{\mathcal{F}}^X}$ is eventually in $\uparrow F$, that is, there exists $(d_0, D_0) \in D_{\mathcal{F}}^X$ such that $ x_{(d,D)}\in \ua F$ for any $(d,D) \ge (d_0, D_0)$.
Hence, there exists a $\mathcal{G}_0\in \mathcal{F}^{(<\omega)}$ such that $\bigcap\{\ua G: G\in \mathcal{G}_0\}\subseteq \ua F$.

(2) $\Rightarrow$ (1): Let $(x_j)_{j \in J} \xrightarrow{GIS} x$. There exists an $\mathcal{F}\subseteq X^{(<\omega)}$ with $\{\ua F : F\in\mathcal F\}\in \ir (P_S(X))$ such that for any $\mathcal{G}\in \mathcal{F}^{(<\omega)}$, $\bigcap\{\ua G: G\in \mathcal{G}\}\neq\emptyset$ and $\ua x\in \cl_{P_S(X,\mathcal{O}_{SI}(X))}\{\ua F: F\in \mathcal{F}\}$, and for any $G \in \mathcal{F}$, $(x_j)_{j \in J}$ is eventually in $\uparrow G$. By hypothesis, there exists an $\mathcal{G}_0\subseteq \mathcal{F}^{(<\omega)}$ such that $\bigcap\{\ua G: G\in \mathcal{G}_0\}\subseteq \ua F$. Since $\mathcal{G}_0$ is finite, we have that $(x_j)_{j \in J}$ is eventually in $\bigcap\{\ua G: G\in \mathcal{G}_0\}$. Hence, $(x_j)_{j \in J}$ is eventually in $\uparrow F$. So we conclude that $F \ll_{GIS} x$.

(2) $\Rightarrow$ (3): Let $\mathcal{F}=\{\ua x: x\in F\}\in \ir (P_S(X))$  which $F$ is an irreducible subset of $X$. Then $\mathcal{F}$ satisfies all the conditions of (2). Hence, there exists an  nonempty finite subset $F_0$ of $F$ such that $\bigcap\{\ua g: g\in F_0\}\subseteq \ua F$ , and thus $F_0^{\ua} \subseteq \ua F$.

\end{proof}

\begin{corollary}
Let $(X, \tau)$ be a $T_0$-space, $x, y \in X$. Then the following statements hold:
\begin{enumerate}[(1)]
    \item If $F \ll_{GIS} G$, then $F \ll G$.
    \item If $X$ is a $GID_C$-space and is a dcpo, then $F \ll G$ implies $F \ll_{GIS} G$.
    In particular, if $P$ is a dcpo, then $F \ll G$ if and only if $F \ll_{GIS} G$ with respect to the Alexandroff topology.
\end{enumerate}
\end{corollary}

\begin{proof}
(1) For any directed set $D$ with $\bigvee D$ exists, $\bigvee D\in \ua G$ implies $\ua G \cap \operatorname{cl}_{\mathcal{O}_{SI(X)}}(D)\neq\emptyset$.   Choose an element $x \in \ua G \cap \operatorname{cl}_{\mathcal{O}_{SI(X)}}(D)$.  It is easy to see that the net $(x_d)_{d \in D}$ with $x_d=d$ for all $d \in D$ $IS$-converges to $x$, and hence also $GIS$-converges to $x$. Hence, $x_d\in \ua F$  holds eventually. Then $D \cap \ua F \neq \emptyset$. Hence, $F \ll G$.

(2) Let $(x_j)_{j \in J} \xrightarrow{GIS} x$ with $x\in \ua G$. There exists an $\mathcal{F}\subseteq X^{(<\omega)}$ with $\{\ua Z : Z\in\mathcal F\}\in \ir (P_S(X))$ such that for any $\mathcal{K}\in \mathcal{F}^{(<\omega)}$, $\bigcap\{\ua K: K\in \mathcal{K}\}\neq\emptyset$ and $\ua x\in \cl_{P_S(X,\mathcal{O}_{SI}(X))}\{\ua Z: Z\in \mathcal{F}\}$, and for any $K \in \mathcal{F}$, $(x_j)_{j \in J}$ is eventually in $\uparrow K$. Since $(X, \tau)$ is  a $GID_C$-space, there exists a filtered subfamily $\mathcal{G}$ with $\mathcal{G}\subseteq \mathbf{Fin}~\!\!X\cap \da_{Q(X)}\{\ua Z : Z\in\mathcal F\}$ such that $\operatorname{cl}_{P_S(X,\mathcal{O}_{SI}(X))}(\mathcal{G}) = \operatorname{cl}_{P_S(X,\mathcal{O}_{SI}(X))}\{\ua Z : Z\in\mathcal F\}$.  Suppose, for the sake of contradiction, that $\ua H\nsubseteq \ua F$ for any $\ua H\in \mathcal{G}$. Let $\mathcal{H}=\{\ua (H\backslash \ua F):\ua H\in \mathcal{G}\}$. Since $\mathcal{G}$ is a filtered subfamily, $\mathcal{H}$ is a filtered subfamily. By Lemma 3.1, there exists a directed set $D\subseteq \bigcup \{H\backslash \ua F:\ua H\in \mathcal{G}\}$ such that $D\cap (H\backslash \ua F)\neq \emptyset$ for all $\ua H\in \mathcal{G}$. Hence, $ \ua x \in\cl_{P_S(X,\mathcal{O}_{SI})}\mathcal{G}=\cl_{P_S(X,\mathcal{O}_{SI})}\{\ua Z : Z\in\mathcal F\}\subseteq\diamond \cl_{\mathcal{O}_{SI}(X)}D=\diamond \da (\vee D)$. Thus, $\vee D\in \ua G$. By $F \ll G$, it follow that $D\cap \ua F\neq\emptyset$ which is in contradiction with $D\subseteq \bigcup \{H\backslash \ua F:\ua H\in \mathcal{G}\}$. Then there exists  $\ua H\in \mathcal{G}$ such that $\ua H\subseteq \ua F$.
\end{proof}

\begin{remark}
(1) The converse of Corollary 5.3(1) may not be true. For example, let $L$ be the $T_0$ space as in Example 4.11. Then $a \ll a$. Since $a$ is in the closure of $\mathbb{N}$ in topological space $SI(L)$, and for any nonempty finite subset $A$ of $\mathbb{N}$, $A^{\ua} \cap \mathbb{N} \neq \emptyset$, we have that $A^{\ua} \nsubseteq \ua a$. So by Theorem 5.2 we conclude that $\{a\} \ll_{GIS} \{a\}$ does not hold.
\end{remark}

\begin{definition}
A $T_0$-space $(X, \tau)$ is called a $GSI$-continuous space if for any $x \in X$, the following conditions hold:
\begin{enumerate}[(1)]
    \item There exists an irreducible set $\mathcal{F}\in \ir P_S(X)$ with $\mathcal{F}\subseteq \{\ua F: F\in w_{GIS}(x)\}$ such that $\ua x \in \operatorname{cl}_{P_S(X,\mathcal{O}_{SI}(X))}\mathcal{F}$;

    \item $\Uparrow_{GIS} F \in \mathcal{O}_{SI}(X)$ for any $F\in w_{GIS}(x)$.
\end{enumerate}
\end{definition}

\begin{proposition}
Let $(X, \tau)$ be a $T_0$-space. Then the following statements are equivalent:

\begin{enumerate}[\rm (1)]
  \item $(X, \tau)$ is a $GSI$-continuous space;
  \item For any net $(x_i)_{i\in I}$ in $X$ and $x \in X$, $(x_i)_{i\in I} \xrightarrow{GIS} x$ if and only if $(x_i)_{i\in I}$ converges to $x$ with respect to the topology $\mathcal{O}_{SI}(X)$.
\end{enumerate}

\end{proposition}

\begin{proof}

(1) $\Rightarrow$ (2): Suppose that  $(x_i)_{i\in I} \xrightarrow{GIS} x$.   As $\mathcal{O}_{SI}(X)\subseteq \mathcal{O}(\mathcal{GIS}(X))$, it follow that  $(x_i)_{i\in I} \xrightarrow{\mathcal{O}_{SI}(X)} x$. Conversely, Suppose that $(x_i)_{i\in I} \xrightarrow{\mathcal{O}_{SI}(X)} x$. Since $X$ is  a $GSI$-continuous space, it follows that $\Uparrow_{GIS} F \in \mathcal{O}_{SI}(X)$ for any $F\in w_{GIS}(x)$ and thus $x_i\in \ua F$ holds eventually.
Moreover, since $X$ is a $GSI$-continuous space, there exists an irreducible set $\mathcal{F}\in \ir P_S(X)$ with $\mathcal{F}\subseteq \{\ua F: F\in w_{GIS}(x)\}$ such that $\ua x \in \operatorname{cl}_{P_S(X,\mathcal{O}_{SI}(X))}\mathcal{F}$. Hence, $\mathcal{F}$ satisfies all the conditions in Definition 4.4. Therefore, $(x_i)_{i\in I} \xrightarrow{GIS} x$.

(2) $\Rightarrow$ (1): Suppose that $x\in X$ and $\mathcal{N}(x)=\{U\in \mathcal{O}_{SI}(X): x\in U\}$.
Let $I=\{(U,e)\in \mathcal{N}(x)\times X: e\in U\}$ and define an order $\leq$ on $I$ as follows:
$(U_1,e_1)\leq (U_2,e_2)$ iff $U_2\subseteq U_1$. Clearly, $I$ is a directed set. For each $(U,e)\in I$, let $x_{(U,e)} = e$.
Clearly, $(x_{(U,e)})_{(U,e)\in I}$ converges to $x$ with respect to the topology $\mathcal{O}_{SI}(X)$.
By condition (2), we have $(x_i)_{i\in I} \xrightarrow{GIS} x$. Then exists $\mathcal{K}\subseteq X^{(<\omega)} $ such that $\{\ua K : K \in \mathcal{K}\} \in \ir(P_S(X))$ and conditions (i)--(iii) in Definition~4.4 are satisfied.
Next, we prove that $F\ll_{GIS}x$ and $\Uparrow_{GIS} F \in \mathcal{O}_{SI}(X)$ for any $F\in \mathcal{K}$.
By $\mathcal{K} \subseteq Elb(x_j)_{j\in J}$, for any $F\in \mathcal{K}$ there is $V\in \mathcal{N}(x)$ such that $x\in V\subseteq \ua F$, and thus $x\in V\subseteq \Uparrow_{GIS} F$ by Remark 4.5(6) and Remark 5.1(5), that is, $F\ll_{GIS}x$. In a similar way, we can prove that $\Uparrow_{GIS} F \in \mathcal{O}_{SI}(X)$.
\end{proof}

From Proposition 4.16 and Proposition 5.6, we obtain the following theorem.

\begin{theorem}
Let $(X, \tau)$ be a $T_0$-space. Then the following statements are equivalent:

\begin{enumerate}[\rm (1)]

  \item $SI(X)$ is a locally hypercompact space;
  \item $(X, \tau)$ is a $GSI$-continuous space;
  \item For any net $(x_i)_{i\in I}$ in $X$ and $x \in X$, $(x_i)_{i\in I} \xrightarrow{GIS} x$ if and only if $(x_i)_{i\in I}$ converges to $x$ with respect to the topology $\mathcal{O}_{SI}(X)$.
\end{enumerate}

\end{theorem}

\begin{definition}
A $T_0$-space $(X, \tau)$ is called a $GIS$-continuous space if for any $x \in X$, the following conditions hold:
\begin{enumerate}[(1)]
    \item There exists an irreducible set $\mathcal{F}\in \ir P_S(X)$ with $\mathcal{F}\subseteq \{\ua F: F\in w_{GIS}(x)\}$ such that $\ua x \in \operatorname{cl}_{P_S(X,\mathcal{O}_{SI}(X))}\mathcal{F}$;

    \item $\Uparrow_{GIS} F \in \mathcal{O}(\mathcal{GIS}(X))$ for any $F\in w_{GIS}(x)$.
\end{enumerate}
\end{definition}

\begin{remark}
Clearly, every $GSI$-continuous space is a $GIS$-continuous space.
\end{remark}

From Theorem 5.7 and Remark 5.9, we obtain the following proposition.

\begin{proposition}
Let $(X, \tau)$ be a $T_0$-space. If $SI(X)$ is a locally hypercompact space, then $X$ is a $GIS$-continuous space.
\end{proposition}

All $GSI$-continuous spaces are $GIS$-continuous, but the converse is not necessarily true. See the following example.

\begin{example}
Let $X$ be an infinite set with the cofinite topology $\tau$. Then $X$ is a $T_1$ space. Since the only irreducible sets admitting a supremum are singletons, it follows that $\tau=\mathcal{O}_{SI}(X)$.
Let $x \in X$, $H\in X^{(<\omega)}$ with $x\in H$ and $\mathcal{F}\subseteq X^{(<\omega)}$ with $\{\ua F : F\in\mathcal F\}\in \ir (P_S(X))$.
And suppose that the point $x$ and the family $\mathcal{F}$ satisfy the conditions in Theorem 5.2(2). As  $\bigcap\{\ua G: G\in \mathcal{G}\}\neq\emptyset$ for any $\mathcal{G}\in \mathcal{F}^{(<\omega)}$, this yields that $\bigcap\{\ua F : F\in\mathcal F\}\neq\emptyset$ and $\bigcap\{\ua F : F\in\mathcal F\}\in X^{(<\omega)}$. Since $\ua x\in \cl_{P_S(X,\mathcal{O}_{SI}(X))}\{\ua F: F\in \mathcal{F}\}$, it follows that $\ua x\supseteq \bigcap\{\ua F: F\in \mathcal{F}\}$ and hence $\{ x\}= \bigcap\{\ua F: F\in \mathcal{F}\}$. As $\mathcal{F}\subseteq X^{(<\omega)}$, there exists a nonempty finite set $\mathcal{G}$ in   $\mathcal{F}$ such that $\bigcap\{\ua G: G\in \mathcal{G}_0\}=\{x\}\subseteq \ua H=H$. By the equivalence of (1) and (2) in Theorem 5.2, we see that $H\ll_{GIS}H$.
Thus, it is easy to know that $\Uparrow_{GIS} H =\ua H\in \mathcal{O}(\mathcal{GIS}(X))$.
Hence, we conclude that $X$ is a $GIS$-continuous space.  However, $\Uparrow_{GIS} H \notin \tau$. Thus, $X$ is not a $GSI$-continuous space.
\end{example}

\begin{proposition}
Let $(X, \tau)$ be an irreducible-complete $GID_C$-space. Then $X$ is an $GIS$-continuous space if and only if for any $x \in X$, the following conditions hold:
\begin{enumerate}[(1)]
    \item $\{\ua G:G\in w_{GIS} (x)\}\in \ir P_S(X)$ and $\ua x =\bigcap \{\ua G:G\in w_{GIS} (x)\}$;
    \item $\Uparrow_{GIS} F \in \mathcal{O}(\mathcal{GIS}(X))$ for any $F\in w_{GIS}(x)$.
\end{enumerate}
\end{proposition}

\begin{proof}
 The sufficiency is obvious by Proposition 3.10.  Let $X$ be a $GIS$-continuous space and $x \in X$. Then there exists an irreducible set  $\mathcal{F}\subseteq \{\ua G: G\in w_{GIS}(x)\}$ such that $\ua x \in \operatorname{cl}_{P_S(X,\mathcal{O}_{SI}(X))}\mathcal{F}$.  It suffices to prove that $\{\ua G: G\in w_{GIS}(x)\}\in \ir P_S(X)$ and $\ua x =\bigcap \{\ua G:G\in w_{GIS} (x)\}$. Since  $x\in \ua G$ for any  $ G\in w_{GIS}(x)$, then $\ua x \subseteq\bigcap \{\ua G:G\in w_{GIS} (x)\}$. Let $W\in \mathcal{O}_{SI}(X)$ and $x\in W$. Thus, $\Box W\cap \mathcal{F}\neq\emptyset$. Therefore, $\bigcap \{\ua G:G\in w_{GIS} (x)\}\subseteq\bigcap\mathcal{F}\subseteq \bigcap\{W\in \mathcal{O}_{SI}(X):x\in W\}=\ua x$.
 Let $U, V \in \tau$ such that $\{\ua G:G\in w_{GIS} (x)\} \cap \Box U \neq \emptyset$ and $\{\ua G:G\in w_{GIS} (x)\} \cap \Box  V \neq \emptyset$.
  Then there exist $G_1, G_2 \in w_{GIS} (x)$ such that $\ua G_1 \subseteq U$, $\ua G_2  \subseteq V$.
  Since $X$ is an $GIS$-continuous space, we have that $\Uparrow_{GIS}  G_1, \Uparrow_{GIS}  G_2 \in \mathcal{O}_{SI}(X)$, and $x \in \Uparrow_{GIS}  G_1$, $x \in \Uparrow_{GIS}  G_2$.
  Then $\ua x \in \operatorname{cl}_{P_S(X,\mathcal{O}_{SI}(X))}\mathcal{F} \cap \Box\Uparrow_{GIS}  G_1$ and $\ua x \in \operatorname{cl}_{P_S(X,\mathcal{O}_{SI}(X))}\mathcal{F} \cap \Box\Uparrow_{GIS}  G_2$.
  By Theorem 4.19, we have that $\mathcal{F} \cap \Box\Uparrow_{GIS}  G_1 \neq \emptyset$ and $\mathcal{F} \cap \Box\Uparrow_{GIS}  G_2 \neq \emptyset$.
  Then $\mathcal{F}  \cap \Box U \neq \emptyset$ and $\mathcal{F} \cap \Box V \neq \emptyset$. Since $\mathcal{F}$ is an irreducible set, we have that $\mathcal{F} \cap\Box  U \cap\Box  V \neq \emptyset$. Hence, $\{\ua G:G\in w_{GIS} (x)\} \cap\Box U \cap\Box V \neq \emptyset$. Therefore, $\{\ua G:G\in w_{GIS} (x)\}$ is an irreducible set.
\end{proof}

\begin{proposition}
Let $(X,\tau)$ be a $GIS$-continuous space. Then the $GIS$-convergence in $X$ is topological with respect to the topology $\mathcal{O}(\mathcal{GIS}(X))$.
\end{proposition}

\begin{proof}
Clearly, the net $(x_j)_{j\in J} \xrightarrow{GIS} x$ implies that $(x_j)_{j\in J}$ converges to $x$ with respect to the topology $\mathcal{O}(\mathcal{GIS}(X))$.
Suppose that $(x_j)_{j\in J}$ converges to $x$ with respect to the topology $\mathcal{O}(\mathcal{GIS}(X))$.
Since $X$ is a $GIS$-continuous space, there exists an irreducible set $\mathcal{F}\in \ir P_S(X)$ with $\mathcal{F}\subseteq \{\ua F: F\in w_{GIS}(x)\}$ such that $\ua x \in \operatorname{cl}_{P_S(X,\mathcal{O}_{SI}(X))}\mathcal{F}$. Let $\ua G\in \mathcal{F}$.
 Then $G \ll_{GIS} x$, and thus $x\in \Uparrow_{GIS} G \in \mathcal{O}(\mathcal{GIS}(X))$.
 Hence $(x_j)_{j\in J}$ is eventually in $\ua G$.
 So we conclude that $(x_j)_{j\in J} \xrightarrow{GIS} x$.
\end{proof}

\begin{proposition}
Let $X$ be a $T_0$ space. If the $GIS$-convergence in $X$ is topological, then $X$ is a $GIS$-continuous space.
\end{proposition}

\begin{proof}
If the $GIS$-convergence in $X$ is topological, then there exists a topology $\tau$ on $X$ such that $(x_j)_{j\in J} \xrightarrow{GIS} x$ if and only if the net $(x_j)_{j\in J}$ converges to $x$ with respect to the topology $\tau$.
Let $x\in X$. Define $I_x=\{(W,a)\mid a\in W,\, x\in W\in \tau\}$. The order $\leq$ on $I_x$ as follows:
$\forall\,(W,a),(V,b)\in I_x,\;(W,a)\leq (V,b)\iff V\subseteq W$.
Then $(I_x,\leq)$ is a preordered set.
Let $(W,a),(V,b)\in I_x$. Then $(W\cap V,x)\in I_x$ and $(W,a),(V,b)\leq (W\cap V,x)$, and thus $I_x$ is a directed set.
Let $x_i = a$ for any $i=(W,a)\in I_x$. Then it is easy to see that $(x_i)_{i\in I_x}$ converges to $x$ with respect to the topology $\tau$, and thus $(x_i)_{i\in I_x}\xrightarrow{GIS} x$.
So we can conclude that there exists an irreducible set $\{\ua F : F\in\mathcal F\subseteq X^{(<\omega)}\}\in \ir (P_S(X))$ such that for any $\mathcal{G}\in \mathcal{F}^{(<\omega)}$, $\bigcap\{\ua G: G\in \mathcal{G}\}\neq\emptyset$ and $\ua x\in \cl_{P_S(X,\mathcal{O}_{SI})}\{\ua F: F\in \mathcal{F}\}$, and for each $\ua F\in \{\ua F: F\in \mathcal{F}\subseteq X^{(<\omega)}\}$, there exists $i_0=(U_e,a)$ such that $x_i\in\ua F$ for all $i=(W,b)\geq i_0$.
Specially, there exists $U_e\in \tau$ such that $x\in U_e\subseteq \ua F$.
Next, we shall prove that $\mathcal{F}\subseteq w_{GIS}(x)$.
Let $F\in \mathcal{F}$, then there exists $U_a\in \tau$ such that $x\in U_a\subseteq \ua F$.
Let $(y_j)_{j\in J}\xrightarrow{GIS} x$. Then $(y_j)_{j\in J}$ converges to $x$ with respect to the topology $\tau$, and thus $(y_j)_{j\in J}$ is eventually in $\ua F$. Hence, $F\ll_{GIS} x$. This shows that $\mathcal{F}\subseteq w_{GIS}(x)$.
 It suffices to prove that $\Uparrow_{GIS}F\in \mathcal{O}(\mathcal{GIS}(X))$. Let $y\in \Uparrow_{GIS}F$. Consider the net $(z_i)_{i\in I_y}$ defined via neighborhoods like the net $(x_i)_{i\in I_x}$ in the beginning of the proof. Then $(z_i)_{i\in I_y}\xrightarrow{GIS} y$.
 It follows from $F\ll_{GIS} y$ that $(z_i)_{i\in I_y}$ is eventually in $\uparrow F$.
 Therefore, there exists $W\in \tau$ such that $y\in W\subseteq \ua F$.
 Let $z\in W$ and $(z_k)_{k\in K}\xrightarrow{GIS} z$.
 Then $(z_k)_{k\in K}$ converges to $z$ with respect to $\tau$, and thus $(z_k)_{k\in K}$ is eventually in $W$. Since $W\subseteq \ua F$, we have that $(z_k)_{k\in K}$ is eventually in $\ua F$.
  Hence, $z\in \Uparrow_{GIS}F$, and so $W\subseteq \Uparrow_{GIS}F$.
  Next, we shall prove that $W\in \mathcal{O}(\mathcal{GIS}(X))$. Let $(p_t)_{t\in T}\xrightarrow{GIS} p\in W$.
   Then $(p_t)_{t\in T}$ converges to $p$ with respect to the topology $\tau$.
    Since $W\in \tau$, we have that $(p_t)_{t\in T}$ is eventually in $W$.
    Then $W\in \mathcal{O}(\mathcal{GIS}(X))$.
     This shows that $\Uparrow_{GIS}F\in \mathcal{O}(\mathcal{GIS}(X))$.
    Therefore, $X$ is a $GIS$-continuous space.
\end{proof}

By Propositions 5.13 and 5.14, we have the following theorem immediately:

\begin{theorem}
Let $(X,\tau)$ be a $T_0$ space. Then the $GIS$-convergence in $X$ is topological if and only if $(X,\tau)$ is a $GIS$-continuous space.
\end{theorem}

\begin{corollary}
Let $P$ be a dcpo. Then the $GS$-convergence in $P$ is topological if and only if $(P,\alpha(P))$ is a $GIS$-continuous space.
\end{corollary}

\vskip 3mm

\noindent{\bf Acknowledgements}

\vskip 3mm

The authors sincerely thank the anonymous referee for his valuable comments and suggestions that have improved this paper substantially.

\end{document}